\documentclass[11pt]{article}
\usepackage[paperwidth=560pt,paperheight=700pt,top=72pt,right=67pt,bottom=72pt,left=72pt]{geometry}
\usepackage{amsmath,amsthm,amsfonts,amssymb,amscd,enumerate,verbatim}
\usepackage[mathscr]{euscript}
\usepackage{array}
\usepackage{hyperref}
\usepackage{fancyhdr}
\usepackage{blindtext}
\usepackage{graphicx,latexsym}
\usepackage{lscape}
\newtheorem{thm}{Theorem}[section]

\newtheorem{lem}[thm]{Lemma}

\theoremstyle{definition}
\newtheorem{definition}{Definition}[section]

\usepackage{latexsym}
\usepackage{amssymb}
\usepackage{amsmath}
\usepackage{graphicx}
\usepackage{hyperref}
\usepackage{caption}
\usepackage{subcaption}
\usepackage{setspace}

\newcommand{\ind}{1\hspace{-2.3mm}{1}}

\renewcommand{\qed}{\hfill{\ \ \rule{2mm}{2mm}} \vspace{0.2in}}

\begin{document}

\date{}
	\title{Minimum Weight Random Graphs with Edge Constraints}
		
	\author{\begin{tabular}{ccl}
		\textbf{Ghurumuruhan Ganesan$^{1}$} \\
		\end{tabular}\\
		\begin{tabular}{c}
			$^{1}$IISER, Bhopal\\
			E-mail: gganesan82@gmail.com
	\end{tabular}}
	\maketitle
\pagestyle{fancy}
\lhead{\textit{Minimum Weight Random Graphs}}
\rhead{\thepage}

\vspace{0.2 cm}
\noindent{\bf Abstract:} In this paper, we study two examples of minimum weight random graphs with edge constraints. First we consider the complete graph on~\(n\) vertices equipped with uniformly heavy edge weights and use iteration methods to obtain deviation estimates for the minimum weight of subtrees with a given number of edges. Next we analyze edge constrained minimum weight paths in the integer lattice~\(\mathbb{Z}^d\) and employ martingale difference techniques to describe the behaviour of the scaled minimum weight in terms of the edge constraint. 
\\[0.5em]
{\bf Keywords:} Minimum spanning tree,  heavy weights, minimum passage time paths, edge constraints. \\[0.5em]
{\bf{2010 Mathematics Subject Classification:}} Primary:
60J10, 60K35; Secondary: 60C05, 62E10, 90B15, 91D30.






\renewcommand{\theequation}{\thesection.\arabic{equation}}
\setcounter{equation}{0}
\section{Introduction} \label{intro}

Trees of complete graphs with random edge weights are important from both theoretical and practical perspectives. For independent and identically distributed (i.i.d.) edge weights with a common cumulative distribution function (cdf)~\(F(.)\) that varies linearly close to zero,~\cite{fre} studied convergence of the weight of the minimum spanning tree (MST) of the complete graph~\(K_n\) on~\(n\) vertices. Later~\cite{ald} studied asymptotics for the \emph{expected} value of the MST weight, when the edge weight distributions follow a power law distribution. The paper~\cite{janson} studied central limit theorems for a scaled and centred version of the MST weight and more recently~\cite{add} studied bounds on the diameter of the MST. The methods involve a combination of graph evolution via Kruskal's agorithm along with a component analysis of random graphs. For MSTs with nonidentical edge weight distributions,~\cite{li} use the Tutte polynomial approach~\cite{ste} to compute expressions for the expected value of~\(MST_n.\)

In Section~\ref{sec_mst} of our paper, we study ``approximate" MSTs containing~\(O(n)\) edges, obtained by placing random heavy weights in each edge of~\(K_n\) that are not necessarily identically distributed but are uniformly heavy. We use stochastic domination to obtain deviation type estimates for the minimum weight and use the martingale method to bound the variance  (see Theorem~\ref{mst_thm}).

Next we study constrained paths in the integer lattice. Consider the following scenario where each edge in the square lattice \(\mathbb{Z}^d\) is associated with a random passage time and it is of interest to determine the minimum passage time~\(T_n\) between the origin and~\((n,0,\ldots,0).\) The case of independent and identically distributed~(i.i.d.) passage times has been well-studied and  detailed results are known regarding the almost sure convergence and convergence in mean of the scaled passage time~\(\frac{T_n}{n}\) (see~\cite{kest1}). Later~\cite{chatterjee} studied central limit theorems for first passage across thin cylinders and recently~\cite{jiang} have studied critical first passage percolation in the triangular lattice.

In many applications, the passage times may not be i.i.d.\ For example, if we model vertices of~\(\mathbb{Z}^d\) as mobiles stations and the edge passage times as the delay in sending a packet between two adjacent stations, then depending on external conditions, the edges may have different passage time distributions. In such cases, it is of interest to study convergence properties of the minimum passage time~\(T_n,\) with appropriate centering and scaling.
In Section~\ref{sec_cons_path} our paper, we state and prove our result  (Theorem~\ref{thm1}) regarding the behaviour of the constrained minimum passage times as a function of the edge constraint.

The paper is organized as follows. In Section~\ref{sec_mst}, we state and prove our result regarding the asymptotic behaviour of weighted trees of the complete graph, with edge constraints. Finally, in Section~\ref{sec_cons_path}, we describe the behaviour of constrained minimum passage time paths in the integer lattice~\(\mathbb{Z}^{d}.\)

\setcounter{equation}{0}
\renewcommand\theequation{\thesection.\arabic{equation}}
\section{Edge Constrained Minimum Weight Trees}\label{sec_mst}
For~\(n \geq 1,\) let~\(K_n \) be the complete graph with vertex set~\(\{1,2,\ldots,n\}.\) Let~\(\{w(i,j)\}_{1 \leq i < j \leq n}\) be independent random variables with corresponding cumulative distribution functions (cdfs)~\(\{F_{i,j}\}_{1 \leq i < j \leq n}\) and for~\(1 \leq j < i \leq n,\) set~\(w(i,j) := w(j,i).\) We define~\(w(e) := w(i,j)\) to be the \emph{weight} of the edge~\(e = (i,j) \in K_n\) and assume throughout that~\(w(e) \leq 1,\) for simplicity. 

A tree in~\(K_n\) is a connected acyclic subgraph. For a tree~\({\cal T}\) with vertex set~\(\{v_1,\ldots,v_t\},\) the weight of~\({\cal T}\) is the sum of the weights of the edges in~\({\cal T};\) i.e.,~\(W({\cal T}) := \sum_{e \in {\cal T}} w(e).\) For~\(1 \leq \tau \leq n-1\) we define
\begin{equation}\label{min_weight_tree}
M_n = M_n(\tau) := \min_{{\cal T}} W({\cal T}),
\end{equation}
where the minimum is taken over all trees~\({\cal T} \subset K_n\) having at least~\(\tau\) edges.

The following is the main result of this section. Constants throughout do not depend on~\(n.\)
\begin{thm}\label{mst_thm} Suppose~\(\tau \geq \rho n\) for some~\(0 < \rho \leq 1\) and also suppose there are positive constants~\(D_1 \leq D_2\) and~\(0 < \alpha < 1\) such that
\begin{equation}\label{dif}
D_1 x^{\frac{1}{\alpha}} \leq F_{i,j}(x) \leq D_2x^{\frac{1}{\alpha}}
\end{equation}
for~\(0 \leq x \leq 1.\) There are positive constants~\(C_i, 1 \leq i \leq 3\) such that
\begin{equation} \label{mst_dev}
\mathbb{P}\left(C_1n^{1-\alpha} \leq M_n \leq C_2 n^{1-\alpha}\right) \geq 1-e^{-C_3n^{1-\alpha}}
\end{equation}
and~\(C_1n^{1-\alpha} \leq \mathbb{E}M_n \leq C_2 n^{1-\alpha}.\) Moreover,~\(var(M_n) \leq 2n.\)
\end{thm}








To prove Theorem~\ref{mst_thm}, we use the following preliminary Lemma regarding the behaviour of the exponential moments of small edge weights. For~\(1 \leq j \leq n\) and distinct deterministic integers~\(1 \leq a_1,\ldots,a_j \leq n\) let
\begin{equation}\label{yj_def}
Y_j = Y_j(a_1,\ldots,a_j) := \min_{a \notin \{a_1,\ldots,a_{j-1}\}} w(a_{j},a).
\end{equation}
\begin{lem}\label{yj_lem} There are positive constants~\(C_1\) and~\(C_2\) not depending on the choice of~\(\{a_i\}\) or~\(j\) such that for any~\(1 \leq j \leq n-1,\)
\begin{equation}\label{eyj_est_lem}
\frac{C_1}{(n-j)^{\alpha}} \leq \mathbb{E}Y_j \leq \frac{C_2}{(n-j)^{\alpha}}.
\end{equation}
Moreover for every~\(s > 1\) there are constants~\(K,C \geq 1\) not depending on the choice of~\(\{a_i\}\) or~\(j\) such that for all~\(1 \leq j \leq n-K\)
\begin{equation}\label{eyj_est_lem2}
\mathbb{E}e^{sY_j} \leq \exp\left(\frac{C}{(n-j)^{\alpha}}\right).
\end{equation}
\end{lem}

\emph{Proof of Lemma~\ref{yj_lem}}: In what follows, we use the following standard deviation estimate. Suppose~\(W_i, 1 \leq i \leq m\) are independent Bernoulli random variables satisfying~\(\mathbb{P}(W_1=1) = 1-\mathbb{P}(W_1~=~0) \leq \mu_2.\) For any~\(0 < \epsilon < \frac{1}{2},\)
\begin{equation}\label{std_dev_up}
\mathbb{P}\left(\sum_{i=1}^{m} W_i > m\mu_2(1+\epsilon) \right) \leq \exp\left(-\frac{\epsilon^2}{4}m\mu_2\right).
\end{equation}
For a proof of~(\ref{std_dev_up}, we refer to Corollary A.1.14, pp. 312 of Alon and Spencer (2008).

We first find the lower bound for~\(\mathbb{E}Y_j\) and then upper bound~\(\mathbb{E}Y_j\) and~\(\mathbb{E}e^{sY_j}\) for constant~\(s > 1\) in that order. The term~\(Y_j\) is the minimum of~\(n-j\) edge weights and so for~\(0 < x < (n-j)^{\alpha}\) we use the upper bound for the cdfs in~(\ref{dif}) to get that~\(\mathbb{P}\left(Y_j > \frac{x}{(n-j)^{\alpha}}\right) \geq \left(1-D_2\frac{x^{\frac{1}{\alpha}}}{n-j}\right)^{n-j},\)
where~\(D_2 \geq 1\) is as in~(\ref{dif}).
Thus
\begin{equation}\label{ey_temp2}
(n-j)^{\alpha} \mathbb{E}Y_j = \int_0^{(n-j)^{\alpha}} \mathbb{P}\left(Y_j > \frac{x}{(n-j)^{\alpha}}\right) dx \geq \int_{0}^{(n-j)^{\alpha}} \left(1-D_2\frac{x^{\frac{1}{\alpha}}}{n-j}\right)^{n-j} dx.
\end{equation}

To evaluate the integral in~(\ref{ey_temp2}), we use~\(1-y \geq e^{-2y}\) for all~\(0 < y < \frac{1}{2}.\) Letting~\(y = \frac{D_2 x^{\frac{1}{\alpha}}}{n-j}\) for~\(0 < x < \left(\frac{n-j}{2D_2}\right)^{\alpha}\) we then have~\((1-y)^{n-j} \geq e^{-2D_2x^{\frac{1}{\alpha}}}\) and substituting this in~(\ref{ey_temp2}) and using~\(D_2 \geq 1,\) we get
\[(n-j)^{\alpha} \mathbb{E}Y_j \geq \int_{0}^{\left(\frac{n-j}{2D_2}\right)^{\alpha}}  e^{-2D_2x^{\frac{1}{\alpha}}} dx \geq \int_{0}^{\left(\frac{1}{2D_2}\right)^{\frac{1}{\alpha}}}e^{-2D_2x^{\frac{1}{\alpha}}} dx  =: C_1\] for all~\(1 \leq  j \leq n-1.\)

For upper bounding~\(\mathbb{E}Y_j,\) we again use the fact that the term~\(Y_j\) is the minimum of~\(n-j\) edge weights and so for~\(0 < x < (n-j)^{\alpha}\) we use the lower bound for the cdfs in~(\ref{dif}) to get~\(\mathbb{P}\left(Y_j > \frac{x}{(n-j)^{\alpha}}\right) \leq \left(1-D_1\frac{x^{\frac{1}{\alpha}}}{n-j}\right)^{n-j} \leq e^{-D_1x^{\frac{1}{\alpha}}}.\)
Thus
\begin{equation}\label{eyj_est}
(n-j)^{\alpha}\mathbb{E}Y_j = \int_0^{(n-j)^{\alpha}} \mathbb{P}\left(Y_j > \frac{x}{(n-j)^{\alpha}}\right) dx \leq \int_0^{\infty} e^{-D_1x^{\frac{1}{\alpha}}} dx  =: C_2,
\end{equation}
a finite positive constant not depending on the choice of~\(\{a_i\}.\)

To compute~\(\mathbb{E}e^{sY_j}\) we split~\(\mathbb{E}e^{sY_j} = I_1 + I_2,\)
where~\(I_1  = \mathbb{E} e^{sY_j} \ind\left(Y_j \leq \frac{1}{2s}\right)\) and~\(I_2 = \mathbb{E}e^{sY_j} \ind\left(Y_j > \frac{1}{2s}\right)\)
and estimate each term separately. To evaluate~\(I_1,\) we bound~\(e^{x} \leq 1+2x\) for~\(x \leq \frac{1}{2}\) and set~\(x= sY_j \leq \frac{1}{2}\) to get that~\(e^{sY_j} \leq 1+2sY_j.\) Thus
\begin{equation}\label{i1_est}
I_1 \leq 1+2s\mathbb{E}Y_j \leq 1+\frac{2sC_2}{(n-j)^{\alpha}},
\end{equation}
using~(\ref{eyj_est}).

To evaluate~\(I_2,\) we recall that~\(Y_j = \min_{a \notin \{a_1,\ldots,a_{j-1}\}} w(a_{j},a) \leq 1\) is the minimum of~\(n-j\) independent edge weights.
Using the upper bound for the cdfs in~(\ref{dif}) we have~\(\mathbb{P}\left(w(a_j,a) > \frac{1}{2s}\right) \leq 1-\frac{D_1}{(2s)^{\frac{1}{\alpha}}} <1,\)
since~\(D_1 \leq 1\) and~\(s~>~1.\) Setting~\(e^{-\theta} = 1-\frac{D_1}{(2s)^{\frac{1}{\alpha}}},\) we therefore have~\(\theta > 0\) and that~\(\mathbb{P}\left(Y_j > \frac{1}{2s}\right) \leq e^{-\theta(n-j)}.\) Finally using~\(Y_j \leq 1\) (since all edge weights are at most one) we get
\begin{equation}\label{i2_est}
I_2 = \mathbb{E}e^{sY_j}\ind\left(Y_j > \frac{1}{2s}\right) \leq e^{s}\mathbb{P}\left(Y_j > \frac{1}{2s}\right) \leq e^{s}e^{-\theta(n-j)} \leq \frac{e^{s}C_2}{(n-j)^{\alpha}},
\end{equation}
where~\(C_2 > 0\) is as in~(\ref{i1_est}), provided~\(n-j \geq K+1\) and~\(K = K(s,C_2)\) is large. The final estimate in~(\ref{i2_est}) is obtained using~\(x^{\alpha} e^{-\theta x} \longrightarrow 0\) as~\(x \rightarrow \infty.\)

From~(\ref{i1_est}) and~(\ref{i2_est}), we therefore get for~\(1 \leq j\leq n-K\) that
\begin{equation}\nonumber
\mathbb{E}e^{sY_j} \leq 1+\frac{2sC_2}{(n-j)^{\alpha}} + \frac{e^{s}C_2}{(n-j)^{\alpha}} \leq \exp\left(\frac{2sC_2+e^{s}C_2}{(n-j)^{\alpha}}\right),
\end{equation}
proving~(\ref{eyj_est_lem2}).~\(\qed\)


\emph{Proof of Theorem~\ref{mst_thm}}: We obtain the lower deviation bound by counting the number of edges with large enough weight and the upper deviation bound by constructing a spanning path with low weight analogous to Aldous (1990). The expectation bounds then follow from~(\ref{mst_dev}). To compute the variance, we use the martingale difference method. In fact, from the variance bound, we get that~\(var\left(\frac{M_n}{n^{1-\alpha}}\right) \leq \frac{C}{n^{1-2\alpha}}\) and so if~\(\alpha < \frac{1}{2},\) then~\(\frac{M_n-\mathbb{E}M_n}{n^{1-\alpha}}\) converges to zero in probability. Details follow.

We begin with the proof of the lower deviation bound. For~\(\gamma > 0\) a small constant,  let~\(R_{tot} := \sum_{e \in K_n} \ind\left(w(e) < \left(\frac{\gamma}{n}\right)^{\alpha}\right)\) be the number of edges of weight at most~\(\left(\frac{\gamma}{n}\right)^{\alpha}.\)
We estimate~\(R_{tot}\) using the standard deviation bound~(\ref{std_dev_up}).
First, we have from the bounds for the cdfs in~(\ref{dif}) that~\(\mathbb{P}\left(w(e) < \left(\frac{\gamma}{n}\right)^{\alpha}\right) <  \frac{D_2\gamma}{n}\) and since the edge weights are independent, we use~(\ref{std_dev_up}) with~\(m = {n \choose 2},\mu_2 = \frac{D_2\gamma}{n}\) and~\(\epsilon  =\frac{1}{4}\)
to get that
\[\mathbb{P}\left(R_{tot} > \frac{5mD_2\gamma}{4n}\right) \leq \exp\left(-m\frac{D_2\gamma}{64n}\right) \leq \exp\left(-\frac{n D_2 \gamma}{256}\right),\]
since~\(m = \frac{n(n-1)}{2} > \frac{n^2}{4}.\) Let~\({\cal T}_n\) be any tree with weight~\(M_n\) and containing at least~\(\tau \geq \rho n\) edges. Using~\(m < \frac{n^2}{2},\) we get that with probability at least\\\(1-e^{-\frac{nD_2\gamma}{256}},\) the weight of~\({\cal T}_n\) is at least
\[\left(\rho n-\frac{5mD_2\gamma}{4n}\right)\cdot \left(\frac{\gamma}{n}\right)^{\alpha} \geq \left(\rho n-\frac{5nD_2\gamma}{8}\right)\cdot \left(\frac{\gamma}{n}\right)^{\alpha}  \geq Cn^{1-\alpha}\] for some constant~\(C > 0,\) provided~\(\gamma >0\) is small. This completes the proof of the lower deviation bound in~(\ref{mst_dev}).

For the upper deviation bound, we consider the spanning path obtained by an incremental approach similar to Aldous (1990). Let~\(i_1 =1\) and among all edges with endvertex~\(i_1,\) let~\(i_2\) be the index such that~\(w(i_1,i_2)\) has the least weight. Similarly, among all edges with endvertex in~\(i_2 \setminus \{i_1\},\) let~\(i_3\) be such that~\(w(i_2,i_3)\) has the least weight. Continuing this way, the path~\({\cal P}_{iter} := (i_1,\ldots,i_n)\)
is a spanning path containing all the nodes and so letting~\(Z_j = w(X_{i_{j-1}},X_{i_j})\)
be the weight of the~\(j^{th}\) edge in~\({\cal P}_{iter},\) we have~\(M_{n} \leq W({\cal P}_{iter}) = \sum_{j=1}^{n} Z_j.\) For~\(s > 0\)
we therefore have
\begin{equation}\label{mst_z}
\mathbb{E}e^{sM_n} \leq \mathbb{E}e^{\sum_{j=1}^{n-1}sZ_j}
\end{equation}
and in what follows, we find an upper bound for the right hand side of~(\ref{mst_z}).

Let~\(a_1 := 1\) and~\(a_l, 2 \leq l \leq j-1\) be deterministic numbers
and suppose the event~\(\{i_1=a_1,\ldots,i_{j-1}=a_{j-1}\}\) occurs
so that~\(Z_l=w(a_l,a_{l+1})\) for~\(1 \leq l \leq j-1\) and~\(Z_j = Y_j = Y_j(a_1,\ldots,a_j) = \min_{a \notin \{a_1,\ldots,a_{j-1}\}} w(a_{j},a)\)
is as in~(\ref{yj_def}). The event~\(\{i_1=a_1,\ldots,i_{j-1}=a_{j-1}\}\) and the random variables~\(w(a_l,a_{l+1}), 1 \leq l \leq j-1\)
depend only the state of edges having at least
one endvertex in~\(\{a_1,\ldots,a_{j-1}\}.\) On the other hand, the random variable~\(Y_j\)
depends only on the state of edges having both endvertices in~\(\{1,\ldots,n\} \setminus \{a_1,\ldots,a_{j-1}\}.\)
Thus
\begin{eqnarray}
&&\mathbb{E}e^{s\sum_{l=1}^{j} Z_l}\ind(i_1=a_1,\ldots,i_{j-1} = a_{j-1}) \nonumber\\
&&\;\;\;\;= \mathbb{E}e^{sY_j} e^{\sum_{l=1}^{j-1} sZ_l} \ind(i_1=a_1,\ldots,i_{j-1}=a_{j-1}) \nonumber\\
&&\;\;\;\;= \mathbb{E}e^{sY_j} \mathbb{E}e^{\sum_{l=1}^{j-1} sZ_l} \ind(i_1=a_1,\ldots,i_{j-1}=a_{j-1}).\label{yj_exp_one}
\end{eqnarray}

Using~(\ref{eyj_est_lem}) we have~\(\mathbb{E}e^{sY_j} \leq \exp\left(\frac{C}{(n-j)^{\alpha}}\right)\)
for all~\(1 \leq j \leq n-K,\) where~\(K\) and~\(C\) do not depend on the choice of~\(\{a_i\}.\)
Thus summing~(\ref{yj_exp_one}) over all possible~\(a_1,\ldots,a_{j-1},\) we get~\(\mathbb{E}e^{s\sum_{l=1}^{j} Z_l} \leq \exp\left(\frac{C}{(n-j)^{\alpha}}\right)\mathbb{E}e^{s\sum_{l=1}^{j-1} Z_l}\)
and continuing iteratively, we get~\(\mathbb{E}e^{s\sum_{l=1}^{j} Z_l} \leq \exp\left(C\sum_{l=1}^{j}\frac{1}{(n-l)^{\alpha}}\right)\)
for~\(n-j \geq K.\) For~\(n-j < K,\) we use the bound~\(\mathbb{E}e^{sY_j} \leq e^{s}\)
since~\(Y_j \leq 1\) (all the edge weights are at most one) and argue as before to get that
\begin{equation}\label{zj_rec3}
\mathbb{E}e^{s\sum_{l=1}^{n-1} Z_l} \leq \exp\left(C\sum_{l=1}^{n-K}\frac{1}{(n-l)^{\alpha}}\right)e^{sK}.
\end{equation}

Comparing with integrals, the term~\[\sum_{l=1}^{n-K}\frac{1}{(n-l)^{\alpha}} =\sum_{j = K+1}^{n-1} \frac{1}{j^{\alpha}} \leq C_3 \int_{K}^{n-1} \frac{1}{x^{\alpha}}dx \leq C_4 n^{1-\alpha}\] for some positive constants~\(C_3,C_4.\) We therefore get from~(\ref{zj_rec3}) and~(\ref{mst_z})
that\\\(\mathbb{E}e^{sM_n} \leq e^{C_5n^{1-\alpha}}\) for some constant~\(C_5 = C_5(s).\) Therefore by Chernoff estimate,
we have~\(\mathbb{P}(M_n \geq C_6n^{1-\alpha}) \leq e^{-C_7 n^{1-\alpha}}\) for some positive constants~\(C_6,C_7.\) This completes the proof of the lower deviation bound in~(\ref{mst_dev}).

Finally, the lower bound on the expectation~\(\mathbb{E}M_n\) follows directly from the lower deviation bound~(\ref{mst_dev}). For the expectation upper bound, we use the fact that the edge weights are at most one and so total weight of any tree containing at least~\(\rho n\) edges is at most~\(n.\) Consequently from the upper deviation bound in~(\ref{mst_dev}), we get that~\(\mathbb{E}M_n \leq C_2 n^{1-\alpha} + n\cdot e^{-Cn^{1-\alpha}} \leq 2C_2n^{1-\alpha}.\) The proof of the variance bound is analogous to the pivotal edge argument in Kesten (1993) together with the fact that the number of edges in a spanning tree is at most~\(n.\) This completes the proof of the Theorem.~\(\qed\)

\setcounter{equation}{0}
\renewcommand\theequation{\thesection.\arabic{equation}}
\section{Edge Constrained Minimum Passage Time Paths} \label{sec_cons_path}
Consider the square lattice \(\mathbb{Z}^d,\) where two vertices~\(w_1 = (w_{1,1},\ldots,w_{1,d})\) and~\(w_2 = (w_{2,1},\ldots,w_{2,d})\) are \emph{adjacent} if~\(\sum_{i=1}^{d}|w_{1,i} - w_{2,i}| = 1\) and adjacent vertices are joined together by an edge. Let~\(\{q_i\}_{i \geq 1}\) denote the set of edges. Each edge~\(q_i\) is equipped with a random passage time \(t(q_i)\) and we define the random sequence~\((t(q_1),t(q_2),\ldots)\) on the probability space~\((\Omega, {\cal F}, \mathbb{P}).\)


A \emph{path}~\(\pi\) is a sequence of distinct adjacent vertices~\((w_1,\ldots,w_{r+1}).\) If~\(e_i, 1 \leq i \leq r\) is the edge with endvertices~\(w_i\)
and~\(w_{i+1},\) then we denote~\(\pi = (e_1,...,e_r).\) By definition,~\(\pi\) is self-avoiding and~\(w_1\) and~\(w_{r+1}\) are said to be the \emph{endvertices} of~\(\pi.\) The length of~\(\pi\) is the number of edges in~\(\pi\) and the passage time of~\(\pi\) is defined as~\(T(\pi) := \sum_{i=1}^{r} t(e_i).\)
\begin{definition}
For~\(k \geq 1\) we define the~\(k-\)\emph{constrained} minimum passage time between the origin and the vertex~\((n,\mathbf{0})\) as~\(T_n(k) := \min_{\pi}T(\pi),\) where the minimum is over all paths~\(\pi\) of length at most~\(k\) and with endvertices~\((0,\mathbf{0})\) and~\((n,\mathbf{0}).\)
We define the \emph{unconstrained} minimum passage time as~\(T_n := \inf_{k \geq 1} T_n(k).\)
\end{definition}
By definition~\(T_n(k) \downarrow T_n\) a.s.\ as~\(k \rightarrow \infty.\) In this section, we are primarily interested in studying how~\(T_n(k)\) varies as the constraint parameter~\(k\) increases and also how fast~\(T_n(k)\) converges to~\(T_n.\)




The following are the main results of this section. Throughout constants do not depend on~\(n.\)
\begin{thm} \label{thm1} Suppose
\begin{equation}\label{sup_i}
\sup_{i \geq 1} \mathbb{P}(t(q_i) \leq \epsilon)  \longrightarrow 0 \text{ as }\epsilon \downarrow 0 \text{ and }\mu_2 := \sup_{i \geq 1} \mathbb{E}t^2(q_i) < \infty.
\end{equation}
\((a)\) There are constants~\(C_1,C_2 > 0\) such that for every~\(k \geq n:\)\\
\begin{equation}\label{bama}
\mathbb{P}\left(C_1n \leq T_n \leq T_n(k) \leq C_2 n\right) \geq 1- \frac{C_2}{n}, \;\;\;\;\;\;\;var(T_n(k)) \leq C_2n
\end{equation}
and~\(C_1n \leq \mathbb{E}T_n \leq \mathbb{E}T_n(k) \leq C_2n.\)\\
\((b)\) There exists a constant~\(C_3 > 0\) such that if~\(k \geq C_3n\) then\\\(\mathbb{P}(T_n \neq T_n(k)) \leq \frac{C_3}{k}.\) If~\(k \geq n^{1+\epsilon}\) for some~\(\epsilon  >0,\) then both~\(\frac{T_n(k)-\mathbb{E}T_n(k)}{n}\) and~\(\frac{T_n-\mathbb{E}T_n}{n}\) converge to zero a.s.\ as~\(n \rightarrow \infty.\)\\
\((c)\) If the edge weights are uniformly square integrable in the sense that\\\(\sup_{i \geq 1} \mathbb{E}t^2(q_i)\ind(t(q_i) \geq M) \longrightarrow 0\) as~\(M \rightarrow \infty,\) then~\(var\left(\frac{T_n}{n}\right) \longrightarrow 0\) as~\(n \rightarrow \infty.\)\\If~\(\sup_{i \geq 1} \mathbb{E}t^p(q_i) < \infty\) for some~\(p > 2,\) then~\(var(T_n) \leq Cn\) for some constant \(C  > 0.\)
\end{thm}

\emph{Proof of Theorem~\ref{thm1}\((a)\)}: Let~\(\mu := \sup_{f} \mathbb{E}t(f)\) be the maximum expected passage time of an edge. We begin by showing that there exists a constant~\(0 < \beta \leq \mu\) such for any integer~\(m \geq 1\) and any path~\(\pi\) containing~\(m\) edges,
\begin{equation}\label{t_pi_ax}
\mathbb{P}\left(T(\pi) \leq \beta m\right) \leq e^{-dm}
\end{equation}
for some positive constant \(\beta  = \beta(d) \leq \mu,\) not depending on~\(m\) or~\(\pi.\) Here~\(d\) is the dimension of the integer lattice under consideration.

Indeed, let~\(\pi = (e_1,\ldots,e_m)\) so that~\(T(\pi) = \sum_{i=1}^{m} t(e_i).\) Using the Chernoff bound we obtain for~\(\delta,s > 0\) that
\begin{equation}
\mathbb{P}(T(\pi) \leq \delta m)  = \mathbb{P}\left(\sum_{i=1}^{m}t(e_i) \leq \delta m\right) \leq e^{s\delta m}\prod_{i=1}^{m}\mathbb{E}\left(e^{-st(e_i)}\right).\label{y_1_eq1}
\end{equation}
For a fixed \(\eta > 0,\) we write~\(\mathbb{E}e^{-st(e_i)} = \int_{t(e_i) < \eta} e^{-st(e_i)} d\mathbb{P} + \int_{t(e_i) \geq \eta} e^{-st(e_i)} d\mathbb{P} \)
and use~\[\int_{t(e_i) < \eta} e^{-st(e_i)} d\mathbb{P} \leq \mathbb{P}(t(e_i) < \eta) \text{ and } \int_{t(e_i) \geq \eta} e^{-st(e_i)} d\mathbb{P} \leq e^{-s\eta} \]
to get that~\[\mathbb{E}e^{-st(e_i)} \leq \mathbb{P}(t(e_i) < \eta) + e^{-s\eta}.\]
Since~\(F(0) = 0,\) we choose~\(\eta > 0\) small so that~\(\mathbb{P}(t(e_i) < \eta) \leq \frac{e^{-6\epsilon}}{2}.\) Fixing such an~\(\eta\) we choose~\(s = s(\eta,\epsilon) > 0\) large so that the second term~\(e^{-s\eta} < \frac{e^{-6\epsilon}}{2}.\) This implies that~\(\mathbb{E}e^{-st(e_i)} \leq e^{-6\epsilon}\) and so from~(\ref{y_1_eq1}) we then get that~\(\mathbb{P}(W({\cal P}) \leq \delta m) \leq e^{s\delta m} e^{-6\epsilon m} \leq e^{-2\epsilon m}\) for all \(m \geq 1,\) provided \(\delta = \delta(s,\epsilon) > 0\) is small. This completes the proof of~(\ref{t_pi_ax}).

Next, for integer~\(m \geq 1\) define the event~\(E_m\) as
\begin{equation}\label{e_k_def}
E_{m} := \bigcap_{r \geq \frac{3\mu}{\beta} m}\;\;\bigcap_{\pi} \left\{T(\pi) \geq \beta r\right\}
\end{equation}
where the second intersection is over all paths with origin as an endvertex and consisting of~\(r\) edges. Thus~\(E_m\) is the event that every path \(\pi\) with origin as an endvertex and consisting of \(r \geq \frac{3 \mu}{\beta}m\) edges has passage time~\(T(\pi) \geq \beta r.\) Since there are at most \((2d)^{r}\) paths of length~\(r\) starting from the origin, the estimate~(\ref{t_pi_ax}) gives
\begin{equation}\label{a_0k}
\mathbb{P}(E_{m}^c)\leq \sum_{r \geq 3\mu\beta^{-1} m} (2d)^{r}e^{-dr} \leq \sum_{r \geq 3\mu\beta^{-1} m} (2e^{-1})^{r} \leq \frac{e^{-\delta m}}{1-2e^{-1}}
\end{equation}
for all \(m \geq 1\) and some positive constant~\(\delta = \delta(d,\mu),\) not depending on~\(m.\) Here, the second inequality in~(\ref{a_0k}) is obtained using the fact that the function~\(xe^{-x} \) attains its maximum at \(x = 1\) and so \(2d e^{-d} \leq 2e^{-1} < 1\) for all \(d \geq 2.\)

Let~\(F_m\) be the event that~\(\sum_{i=1}^{m}t(f_i) \leq 2\mu m\) where~\(f_i\) is the horizontal edge with endvertices~\((i-1,\mathbf{0})\) and \((i,\mathbf{0}).\) Letting~\(X_i := t(f_i) - \mathbb{E}t(f_i)\) and using the fact that~\(\{X_i\}\) are independent, we then get from Chebychev's inequality that
\begin{equation}\label{fm_est}
\mathbb{P}(F_m^c) \leq \mathbb{P}\left(\sum_{i=1}^{m}X_i \geq \mu m\right) \leq \frac{\sum_{i=1}^{m} var(X_i)}{\mu^2m^2} \leq \frac{C}{m}
\end{equation}
for some constant~\(C > 0.\) Now set~\(m = \frac{\beta n}{3\mu} <n\) and suppose~\(E_m \cap F_n\) occurs. From~(\ref{a_0k}),~(\ref{fm_est}) and the union bound, we get that~\(\mathbb{P}(E_m \cap F_n) \geq 1-\frac{C}{n}\) for some constant~\(C > 0.\) Since~\(F_n\) occurs, we get that~\(T_n(k) \leq 2\mu n\) and since~\(E_m\) occurs, we get that any path starting from the origin and containing~\(r \geq \frac{3\mu}{\beta}m = n\) edges has weight at least~\(\beta r \geq 3\mu m = \beta n.\) This obtains the first estimate~(\ref{bama}).

Next, using the bounded second moment assumption in~(\ref{sup_i}) and arguing as in the variance estimate in Theorem~\(1,\) Kesten (1993) we get that~\(var(T_n(k)) \leq C \mathbb{E}N_n(k),\) where~\(N_n(k)\) is the number of edges in the path with passage time~\(T_n(k).\) If~\(E_m \cap F_n\) occurs, then from the discussion in the previous paragraph, we get that~\(N_n(k) \leq \frac{3\mu}{\beta}n.\)
For~\(x \geq \frac{3\mu}{\beta}n\) we assume for simplicity that~\(y = \frac{\beta x}{3\mu}\) is an integer and write
\begin{eqnarray}
\mathbb{P}(N_n(k)  \geq x) &\leq& \mathbb{P}\left(\{N_n(k) \geq x\} \cap E_y\right) + \mathbb{P}\left(E_y^c\right) \nonumber\\
&\leq& \mathbb{P}\left(\{N_n(k) \geq x\} \cap E_m\right) + D_1 e^{-D_2 x}\label{thmoo_one}
\end{eqnarray}
for some constants~\(D_1,D_2> 0 \) by~(\ref{a_0k}). If~\(N_n(k) \geq x\) and the event~\(E_y\) occurs, then every path containing~\(r \geq \frac{3\mu}{\beta}y = x\) edges has weight at least~\(\beta r \geq \beta x.\) Thus~\(T_n(k) \geq \beta x\) and so
\[\mathbb{P}\left(\{N_n(k)  \geq x\} \cap E_y\right) \leq \mathbb{P}(T_n(k) \geq \beta x) \leq \mathbb{P}\left(\sum_{i=1}^{n}t(f_i) \geq \beta x\right).\]

Since~\(\beta x \geq 3\mu n\) we have that~\(\beta x - \sum_{i=1}^{n} \mathbb{E}t(f_i) \geq \beta x - \mu n \geq \frac{2\beta x}{3}.\)
Consequently, recalling that~\(X_i = t(f_i) -\mathbb{E}t(f_i)\) and using the fact that~\(var(X_i) \leq \mathbb{E}t^2(f_i) \leq C\) for some constant~\(C > 0,\) we get from Chebychev's inequality that
\begin{equation}
\mathbb{P}\left(\{N_n(k)  \geq x\} \cap E_y\right) \leq \mathbb{P}\left(\sum_{i=1}^{n}X_i \geq \frac{2\beta x}{3}\right) \leq \frac{ D_3\sum_{i=1}^{n} var(X_i)}{x^2} \leq \frac{D_4n}{x^2} \label{thmoo_two}
\end{equation}
for some constants~\(D_3,D_4 > 0.\) Combining~(\ref{thmoo_one}) and~(\ref{thmoo_two}), we get that~\(\mathbb{P}(N_n(k) \geq x) \leq \frac{D_5}{x^2}\) for~\(x \geq \frac{3\mu}{\beta}n\) and so~\(\mathbb{E}N_n(k) \leq D_6 n\) for some constant~\(D_6 > 0.\) Plugging this into the variance estimate for~\(T_n(k)\) obtained in the previous paragraph, we get the second estimate in~(\ref{bama}).

The lower expectation bound follows directly from the lower deviation bound in~(\ref{bama}). The upper expectation bound follows from the fact that~\(\mathbb{E}T_n \leq \sum_{i=1}^{n} \mathbb{E}t(f_i) \leq \mu n.\) This completes the proof of part~\((a).\)~\(\qed\)

\emph{Proof of Theorem~\ref{thm1}\((b)\)}: Let~\(\beta,\mu\) be as in part~\((a)\) and for~\(k \geq n\) suppose that the event~\(E_k \cap F_k\) occurs. From the discussion following~(\ref{fm_est}), we know that~\(\mathbb{P}(E_k \cap F_k) \geq 1-\frac{C}{k}\) for some constant~\(C >0.\) The minimum passage time between the origin and~\((n,\mathbf{0})\) is at most~\(2\mu k\) and any path starting from the origin and containing~\(r \geq  \frac{3 \mu}{\beta}k\) edges has passage time at least~\(\beta r \geq 3\mu k.\) Thus~\(T_n\left(\frac{3\mu k}{\beta}\right) = T_n\) and this obtains the probability estimate in~\((b)\) with~\(C_3~=~\frac{3\mu}{\beta}.\)

We prove the a.s.\ convergence in two steps. In the first step, we use a subsequence argument to show that~\(\frac{T_n(k)-\mathbb{E}T_n(k)}{n}\) converges to zero a.s. In the second step we show that the \emph{difference}~\(\frac{T_n(k)-T_n}{n}\) converges to zero a.s.\ and in~\(L^1,\) provided~\(k\) is sufficiently large. This then obtains the a.s.\ convergence for~\(\frac{T_n-\mathbb{E}T_n}{n}.\)

\underline{\emph{Step 1}}: We begin with a description of the sub-additivity property of the unconstrained passage time~\(T_n.\) If~\(T_{n,m}\) is the minimum passage time between~\((n,\mathbf{0})\) and~\((m,\mathbf{0}),\) then~\(T_n \leq T_m + T_{n,m}.\) This is because the concatenation of the minimum passage time path with endvertices~\((0,\mathbf{0})\) and~\((m,\mathbf{0})\) and the minimum passage time path with endvertices~\((m,\mathbf{0})\) and~\((n,\mathbf{0})\) contains a path with endvertices~\((0,\mathbf{0})\) and~\((n,\mathbf{0}).\) Switching the roles of~\(m\) and~\(n\) we therefore have that~\(|T_n-T_m| \leq T_{n,m}\) and we refer to this estimate as the \emph{sub-additive property} of~\(T_n.\)

Letting~\(k \geq n^{1+\epsilon}\) with~\(\epsilon > 0,\) we now perform the subsequence argument.  Setting~\(U_n := T_n(k)\) and~\(S_n := U_n - \mathbb{E}U_n,\) we first show that~\(\frac{S_n}{n} \longrightarrow 0\) a.s.\ as~\(n~\rightarrow~\infty.\) Indeed, from the variance estimate for~\(U_n\) in~(\ref{bama}), we know that~\(\mathbb{E}S^2_n \leq C n\) for some constant \(C > 0\) and so for a fixed \(\delta  > 0,\) the sum \[\sum_{n \geq 1} \mathbb{P}(|S_{n^2}| > n^2 \delta) \leq \sum_{n \geq 1}\frac{\mathbb{E}S^2_{n^2}}{\delta^2 n^4} \leq \sum_{n \geq 1} \frac{C}{\delta^2 n^{2}} < \infty.\] Since this is true for all \(\delta > 0,\) Borel-Cantelli Lemma implies that~\(\frac{S_{n^2}}{n^2} \) converges to zero a.s.\ as \(n \rightarrow \infty.\)

To estimate the intermediate values of~\(S_j,\) we let~\(n^2 \leq j < (n+1)^2\) and set~\(R_n := \max_{n^2 \leq j < (n+1)^2} |S_j-S_{n^2}|\) and show below that~\(\frac{R_n}{n^2} \longrightarrow 0\) a.s.\ as~\(n~\rightarrow~\infty.\)  This would imply that for \(n^2 \leq j < (n+1)^2, \) \[\frac{|S_j|}{j}  \leq \frac{|S_j-S_{n^2}|}{j} + \frac{|S_{n^2}|}{j} \leq \frac{|S_j-S_{n^2}|}{n^2} + \frac{|S_{n^2}|}{n^2} \leq \frac{D_{n}}{n^2} + \frac{|S_{n^2}|}{n^2}\] and so~\(\frac{S_j}{j}\) converges to zero a.s.\ as \(j \rightarrow \infty.\)

To estimate~\(R_n,\) use first the triangle inequality to get that
\begin{equation}
|S_j - S_{n^2}| \leq |U_j- U_{n^2}| + \mathbb{E}|U_j - U_{n^2}| \label{eq_s2}
\end{equation}
We know that~\(\mathbb{P}(T_j \neq U_j) \leq \frac{D_1}{k(j)} \leq \frac{D_2}{j^{1+\epsilon}}\) for some constants~\(D_1,D_2 > 0\) and so setting~\(E_{tot} := \bigcap_{j=n^2}^{(n+1)^2}\{T_j = U_j\},\) we get from the union bound that
\begin{equation}\label{e_tot_est}
\mathbb{P}(E_{tot}) \geq 1- \sum_{j=n^2}^{(n+1)^2} \frac{D_2}{j^{1+\epsilon}} \geq 1-\frac{D_3 n}{n^{2+2\epsilon}} = 1-\frac{D_3}{n^{1+2\epsilon}}
\end{equation}
for some constant~\(D_3 > 0.\)

We now write~\(|U_j-U_{n^2}| = |T_j-T_{n^2}| \ind(E_{tot}) + |U_j-U_{n^2}| \ind(E^c_{tot})\) and evaluate each term separately.
For~\(n^2 \leq j < (n+1)^2\) we know by the subadditivity property that~\(|T_j-T_{n^2}| \leq T_{j,n^2} \leq \sum_{i=n^2}^{(n+1)^2}t(f_i) =: A_n\)
and by definition, we have that~\(U_j \leq \sum_{i=1}^{(n+1)^2} t(f_i) =: J_n.\) Thus~\(|U_j-U_{n^2}| \leq A_n + 2J_n\ind(E^c_{tot})\) for all~\(n^2 \leq j \leq(n+1)^2\) and so from~(\ref{eq_s2}), we see that
\begin{equation}\label{dn_est}
R_n \leq A_n + 2J_n\ind(E^c_{tot}) + \mathbb{E}A_n + 2\mathbb{E}J_n \ind(E^c_{tot}).
\end{equation}

Based on~(\ref{dn_est}), it suffices to show that both the terms~\(\frac{A_n}{n^2}\) and~\(\frac{J_n \ind(E^c_{tot})}{n^2}\) converge to zero a.s.\ and in~\(L^1\) as~\(n \rightarrow \infty.\) First, from the estimate~\(\mathbb{E}A_n \leq \sum_{i=n^2}^{(n+1)^2}\mathbb{E}t(f_i) \leq Cn\) for some constant~\(C >0, \) we get that~\(\frac{\mathbb{E}A_n}{n^2} \longrightarrow 0\) as~\(n~\rightarrow~\infty.\) Next, let~\(0 < \theta < 1\) be any constant. Using Chebychev's inequality and arguing as in~(\ref{fm_est}), we get that~\(\mathbb{P}(A_n \geq 2\mu n^{1+\theta}) \leq \frac{C}{n^{1+2\theta}}\) for all~\(n\) large and so by the Borel-Cantelli Lemma we get that~\(\mathbb{P}\left(A_n \leq 2\mu n^{1+\theta} \text{ for all large }n\right) = 1.\) Since~\(\theta < 1,\) we get that~\(\frac{A_n}{n^2} \longrightarrow 0\) a.s.\ as~\(n \rightarrow \infty.\)

To evaluate~\(J_n\ind(E^c_{tot}),\) we use~(\ref{e_tot_est}) and the Borel-Cantelli Lemma to get~\(\ind(E^c_{tot}) \longrightarrow 0\) a.s.\ as~\(n \rightarrow \infty.\) Thus~\(\frac{J_n\ind(E^c_{tot})}{n^2} \longrightarrow 0\) a.s.\ as~\(n \rightarrow \infty.\)
Using the Cauchy-Schwartz inequality, we also get that~\(\mathbb{E}J_n\ind(E^c_{tot}) \leq \left(\mathbb{E}J_n^2\right)^{1/2} \left(\mathbb{P}(E^c_{tot})\right)^{1/2}.\) By the AM-GM inequality and the bounded second moment assumption in~(\ref{sup_i}), we get that~\[\mathbb{E}J_n^2 \leq (n+1)^2 \sum_{i=1}^{(n+1)^2} \mathbb{E}t^2(f_i) \leq C_2 n^4\] for some constant~\(C_2 > 0\) and so using~(\ref{e_tot_est}), we get that~\[\mathbb{E}J_n\ind(E^c_{tot})  \leq \sqrt{C_2} n^2 \cdot \left(\frac{D}{n^{1+2\epsilon}}\right)^{1/2}.\] Thus~\(\frac{\mathbb{E}J_n\ind(E^c_{tot})}{n^2}~\longrightarrow~0\) as~\(n \rightarrow \infty\) as well and this completes the proof of~\(\frac{S_n}{n} \longrightarrow 0\) a.s.\ as~\(n \rightarrow \infty.\)

\underline{\emph{Step 2}}: We now set~\(k = n^{1+\epsilon}\) and show that~\(\frac{T_n(k)-T_n}{n}\) converges to zero a.s.\ and in~\(L^1.\) First using~\(\mathbb{P}(T_n(k) \neq T_n) \leq \frac{D_1}{n^{1+\epsilon}}\) for some constant~\(D_1 > 0,\) we get from the Borel-Cantelli Lemma that~\(\frac{T_n(k)-T_n}{n} \longrightarrow 0\) a.s.\ as~\(n \rightarrow \infty.\) Next using the fact that~\(T_n(k)\) and~\(T_n\) are both bounded above by~\(\sum_{i=1}^{n} t(f_i)\) we have that~\(\mathbb{E}|T_n(k)-T_n| = \mathbb{E}|T_n(k)-T_n|\ind(T_n(k) \neq T_n)\) is bounded above by
\begin{eqnarray}
\mathbb{E}\sum_{i=1}^{n}t(f_i) \ind(T_n(k) \neq T_n) &\leq& \mathbb{E}^{1/2} \left(\sum_{i=1}^{n} t(f_i)\right)^2 \mathbb{P}^{1/2}(T_n(k) \neq T_n) \nonumber\\
&\leq& \mathbb{E}^{1/2} \left(\sum_{i=1}^{n} t(f_i)\right)^2 \left(\frac{D_1}{n^{1+\epsilon}}\right)^{1/2}. \nonumber
\end{eqnarray}
Using~\((\sum_{i=1}^{l}a_i)^2 \leq \l \sum_{i}a_i^2\) we have~\(\mathbb{E}\left(\sum_{i=1}^{n} t(f_i)\right)^2 \leq n \sum_{i=1}^{n} \mathbb{E}t^2(f_i) \leq D_2n^2\) for some constant~\(D_2 > 0.\) Thus~\(\mathbb{E}|T_n(k)-T_n| \leq D_3 (n^{1-\epsilon})^{1/2} = o(n)\) and so\\\(\frac{\mathbb{E}|T_n(k)-T_n|}{n} \longrightarrow 0\) as~\(n \rightarrow \infty.\) This completes the proof of a.s.\ convergence in part~\((b).\)~\(\qed\)



\emph{Proof of Theorem~\ref{thm1}\((c)\)}: Using~\((a+b)^2 \leq 2(a^2 + b^2)\) for any two real numbers~\(a\) and~\(b\) we have that the variance of the sum of any two random variables~\(X\) and~\(Y\) satisfies
\begin{equation}\label{var_xy}
var (X+Y) = \mathbb{E}\left( (X - \mathbb{E}X) + (Y - \mathbb{E}Y)\right)^2 \leq 2(var(X) + var(Y)).
\end{equation}
Setting~\(T = T_n, U = T_n(k), X = T-U\) and \(Y = {U}\) we get that
\begin{equation}\label{var_t}
var(T) \leq 2 var(T - U)  + 2 var(U) \leq 2 var(T-U) + D_1 n
\end{equation}
for all \(n\) large and some constant~\(D_1 >0,\) by the variance estimate for~\(U\) in part~\((a)\) of this Theorem.

We estimate~\(var(T-U)\) as follows. Using~\(T \leq U\) we write
\[\mathbb{E}(T-U)^2 = \mathbb{E}(T-U)^2 \ind(T \neq U) \leq 2\mathbb{E}(T^2+U^2)\ind(T \neq U) \leq 4\mathbb{E}U^2 \ind(T \neq U).\]
Since~\(U \leq \sum_{i=1}^{n}t(f_i),\) we have that~\(U^2 \leq n \sum_{i=1}^{n}t^2(f_i)\) and so~\(var(T-U)\) is bounded above by
\begin{equation}\label{etu}
\mathbb{E}(T-U)^2 \leq 4n\sum_{i=1}^{n}\mathbb{E}t^2(f_i)\ind(T \neq U) \leq 4n^2 \sup_{i} \mathbb{E}t^2(f_i) \ind(T \neq U).
\end{equation}
Let~\(\theta > 0\) be a constant and split
\begin{eqnarray}\label{eq_zn3}
\mathbb{E}t^2(f_i)\ind (T \neq U) &=& \mathbb{E}t^2(f_i)\ind(\{T \neq U\} \cap \{t(f_i) < n^{\theta}\}) \nonumber\\
&&\;\;\;\;\;\;\;+\;\;\;\mathbb{E}t^2(f_i)\ind(\{T \neq U\} \cap \{t(f_i) \geq n^{\theta}\}).
\end{eqnarray}
From~(\ref{bama}), we know that there are constants~\(D_2,D_3 > 0\) such that if~\(k \geq D_2n\) then\\\(\mathbb{P}(T \neq U) \leq \frac{D_3}{k}.\) With this choice of~\(k,\) the first term in~(\ref{eq_zn3}) is bounded above by
\begin{equation}\label{ondra}\mathbb{E}t^2(f_i)\ind(\{T \neq U\} \cap \{t(f_i) < n^{\theta}\}) \leq n^{2\theta}\mathbb{P}(T \neq U) \leq \frac{D_2n^{2\theta}}{k} \leq \frac{1}{n^3}
\end{equation}
provided we choose~\(k\) larger if necessary so that~\(k \geq n^{2\theta}\log{n}.\)

We now consider the case where the uniform square integrability condition holds. For any~\(\eta > 0\) and all~\(n\) large, the final term in~(\ref{eq_zn3}) is then at most~\(\mathbb{E}t^2(f_i) \ind(t(f_i) \geq n^{\theta}) \leq \eta.\) Combining this estimate with~(\ref{ondra}), we get that~\(\mathbb{E}t^2(f_i)\ind (T \neq U) \leq \frac{1}{n^3} + \eta \leq 2\eta\) for all~\(n\) large and so from~(\ref{etu}) we get that~\(\mathbb{E}(T-U)^2 \leq 8n^2 \eta.\) Plugging this into~(\ref{var_t}) we get that~\(var(T) \leq D_1 n  + 8n^2 \eta \) and since \(\eta > 0\) is arbitrary, this implies that~\(var\left(\frac{T}{n}\right) = o(1).\)

Suppose now that bounded~\(p^{th}\) moment condition holds for some~\(p > 2.\) Using H\"older's inequality and Markov inequality in succession, the final term in~(\ref{eq_zn3}) is at most
\begin{eqnarray}
\mathbb{E}t^2(f_i) \ind(t(f_i) \geq n^{\theta}) &\leq& \left(\mathbb{E}t^{p}(f_i)\right)^{2/p} \mathbb{P}\left(t(f_i) \geq n^{\theta}\right)^{1-2/p} \nonumber\\
&\leq& \left(\mathbb{E}t^{p}(f_i)\right)^{2/p} \left(\frac{\mathbb{E}t^{p}(f_i)}{n^{\theta p}}\right)^{1-2/p} \nonumber\\
&\leq& \frac{D_4}{n^{\theta(p-2)}} \label{renda}
\end{eqnarray}
for some constant~\(D_4 > 0,\) by the bounded~\(p^{th}\) moment assumption. We choose~\(\theta > 0\) large so that the final term in~(\ref{renda}) is at most~\(\frac{1}{n^3}.\) With this choice of~\(\theta\) we get from~(\ref{ondra}) that~\(\mathbb{E}t^2(f_i)\ind (T \neq U) \leq \frac{2}{n^3}\) and plugging this into~(\ref{etu}), we get that~\(\mathbb{E}(T-U)^2 \leq \frac{D_5}{n}\) for constant~\(D_5 >0.\)  From~(\ref{var_t}), we then get that~\(var(T) \leq D_6n\) for some constant~\(D_6 > 0.\) This completes the proof of part~\((c).\)~\(\qed\)

\underline{\emph{Remark}}: For~\(p > 2,\) the bounded~\(p^{th}\) moment condition in Theorem~\ref{thm1}\((c)\) is stronger than the uniformly square integrable condition which in turn is stronger than the bounded second moment condition in~(\ref{sup_i}). For the particular case of i.i.d.\ passage times, uniform square integrability is implied by the bounded second moment condition and the first condition in~(\ref{sup_i}) above simply states that the passage times are a.s.\ positive.

\subsection*{Acknowledgement}
I thank Professors Rahul Roy, C. R. Subramanian and the referee for crucial comments that led to an improvement of the paper. I also thank IMSc and IISER Bhopal for my fellowships.





\bibliographystyle{plain}

\end{document}